\documentclass[a4paper,12pt]{article}
\usepackage{latexsym}
\usepackage{amsmath}
\usepackage{amssymb}
\usepackage{graphicx}
\usepackage[round]{natbib}
\usepackage{delarray}
\usepackage[usenames]{color}
\usepackage{colortbl}
\usepackage{theorem}
\usepackage{mathrsfs}
{\theorembodyfont{\upshape}\newtheorem{rem}{Remark}[section]}
\newtheorem{cor}{Corollary}[section]
\newtheorem{df}{Definition}[section]
\newtheorem{lem}{Lemma}[section]
{\theorembodyfont{\upshape}
}
 
\newtheorem{thm}{Theorem}[section]
\newtheorem{prop}{Proposition}[section]
{\theorembodyfont{\upshape}\newtheorem{ex}{Example}[section]}
\def\rnum#1{\expandafter{\romannumeral #1}} 
\def\Rnum#1{\uppercase\expandafter{\romannumeral #1}}

\title{Generalized fractional Ornstein-Uhlenbeck processes}
\author{Kotaro ENDO${}^{\dag}$ and Muneya MATSUI
\footnote{Correspondence to Department of Mathematics, Keio University,
3-14-1 Hiyoshi Kohoku-ku, Yokohama 223-8522, Japan.
 E-mail:mmuneya@gmail.com.\ 
This research is supported by JSPS Research Fellowships for Young Scientists.}
\\
$\dag$ 
\small Department of Mathematics, Keio University \\
$\ddag$
\small Center for Mathematical Science, Munich University of Technology \\
\small and Department of Mathematics, Keio University}
\begin{document}
\maketitle
 
\begin{abstract}
 We introduce an extended version of the fractional Ornstein-Uhlenbeck (FOU)
 process where the integrand is replaced by the exponential of an independent L\'evy process.  
 We call the process the generalized fractional Ornstein-Uhlenbeck (GFOU)
  process. Alternatively, the process can be constructed 
 from a generalized Ornstein-Uhlenbeck (GOU) process using an
 independent fractional Brownian motion (FBM) as integrator.
 We show that the GFOU process is well-defined by checking the existence of the
 integral included in the process, and investigate its properties.
 It is proved that the process has a stationary version and exhibits long
 memory.  We also find that the process satisfies
 a certain stochastic differential equation. Our underlying intention is
 to introduce long memory into the GOU process which has short
 memory without losing the possibility of jumps. 
 Note that both FOU and GOU processes have found application in
 a variety of fields as useful alternatives to the Ornstein-Uhlenbeck (OU)
  process.
\end{abstract}

\section{Introduction}
The fractional Brownian motion (FBM) is one of the most popular processes for
constructing long-range dependent stochastic processes with continuous path 
and its fields of applications are very wide. To name just a few, we see FBM
models in the fields of telecommunications, signal processes,
environmental models and economics.  
A recent reference is e.g., \cite{Doukhan-etal:2003}.
Statistical methods for FBM 
have also been studied (see e.g. \cite{Beran:1994}). 

We review the definition and name properties of FBM.
\begin{df}
Let $0<H\le1$. A fractional Brownian motion $B^H:=\{B_t^H\}_{t\in \mathbb{R}}$
 is a centered Gaussian process with $B_0^H=0$ and 
$\mathrm{Cov}\left(B_t^H,B_s^H\right)=\frac{1}{2}\left(|t|^{2H}+|s|^{2H}-|t-s|^{2H}\right),\
 t,s\in\mathbb{
R}.$
\end{df}
From the definition FBM  has stationary
 increments and is self-similar with index $H$, i.e., for
 $c>0$
 $\{B_{ct}^H\}_{t\in\mathbb{R}}\stackrel{d}{=}\{c^HB_t^H\}_{t\in\mathbb{R}}$
 where $\stackrel{d}{=}$ denotes equality of all finite
 dimensional distributions.
 While $\{B_t^H\}_{t\in\mathbb{R}}$ with $H=1/2$ is a two-sided Brownian
 motion (BM) and has
 independent increments, $\{B_t^H\}_{t\in\mathbb{R}}$ with
 $H\in(0,\frac{1}{2})\cup(\frac{1}{2},1]$ has dependent
 increments. For $0<h<s$ and $t \in\mathbb{R}$ and $N\in \mathbb N$,
{\allowdisplaybreaks
\begin{eqnarray*}
\Gamma_h(s)&:=&\text{Cov}\left(B_{t+h}-B_t,B_{t+s+h}-B_{t+s}\right)
=\text{Cov}\left(B_{h},B_{s+h}-B_{s}\right)\\
&=& \sum_{n=1}^\infty
 \left(\prod_{k=0}^{2n-1}\frac{h^{2n}}{(2n)!}\left(2H-k\right)\right) 
s^{2H-2n}\\
&=& \sum_{n=1}^N
 \left(\prod_{k=0}^{2n-1}\frac{h^{2n}}{(2n)!}\left(2H-k\right)\right) 
s^{2H-2n}+O\left(s^{2H-2N-2}\right),\ \text{as}\ s\to\infty.
\end{eqnarray*}}
Thus, $B^H$ with $H\in(\frac{1}{2},1]$ has a long
memory property, namely,
$
\sum_{n=0}^\infty \Gamma_h(nh)=\infty.
$ 
Finally $B^H$ is known to have bounded $p$-variation for $1/H<p<\infty$
(see Proposition 2.2 of \cite{Mikosch:Norvaisa:2000}). 
For a more detailed theoretical treatment, we refer to \cite{embrecht:maejima:2002}
or  \cite{Samorodnitsky:Taqqu:1994}. 

On the other hand, recently, extensions of
the classical OU process have been suggested mainly on demand of
applications. The generalized Ornstein-Uhlenbeck process given below is one with L\'evy processes and plays
important role in Economics (pricing of Asian options, perpetuities and risk theory).   
For its theories and applications, we refer to e.g, 
\cite{carmona:petit;yor:2001}, \cite{erickson:maller:2004},
\cite{lindner:maller:2005} and \cite{Kuluppelberg:Kostadinova:2008}.
A multivariate extension is also considered in
\cite{kondo:maejima:sato:2006} and \cite{endo:matsui:2008}.
Let $\{(\xi_t,\eta_t)\}_{t \ge 0}$ be a bivariate L\'evy process and
$V_0$ be an independent initial random variable. 
A generalized Ornstein-Uhlenbeck (GOU)
process is defined as
\begin{equation}
\label{eq:def:GOU}
V_t=e^{-\xi_t}\left(V_0+\int_0^te^{\xi_{s-}}d\eta_s\right),\ t\ge0.
\end{equation}
The stationarity and the convergence or divergence property have been
intensively studied.
If $\{\xi_t\}_{t \ge 0}$ and $\{\eta_t\}_{t \ge 0}$ are independent and $\int_0^t
e^{-\xi_{s-}}d\eta_s$ converges $a.s.$ as $t\to\infty$ to a finite random variable,
$V:=\{V_t\}_{t\ge0}$ has the stationary version (see e.g. Remark 2.2 of \cite{lindner:maller:2005}). 
The short memory property of $V$ was also shown in Section 4 of \cite{lindner:maller:2005}.

Another extension of the original OU process is the fractional
Ornstein-Uhlenbeck process, where FBM is used as integrator. An advantage
of using the process is to realize stationary long range dependent processes. 
Let $\lambda>0$ and an initial random variable $X_0^H\in L^1$. 
A fractional Ornstein-Uhlenbeck (FOU) process is defined as 
\begin{equation}
\label{eq:def:FOU}
X_t^{H}=e^{-\lambda t}\left(X^H_0+ \int_0^t e^{\lambda s}dB_s^H\right).
\end{equation}
Here we need a non-semimartingale approach to construct
stochastic integrals with FBM. We can find several useful theoretical
tools in e.g., \cite{lin:1995}, \cite{Mikosch:Norvaisa:2000} or
\cite{pipiras:taqq:2000b}.  
\cite{cheridito:kawaguchi:maejima:2003} has shown that the FOU process is the unique continuous solution of 
a Langevin equation: 
$X_t=X_0^H-\lambda\int_{0}^{t}X_sds+B_t^H,\
t\ge0$ and 
investigated its dependence properties.
The main purpose of this paper is to construct a version of the
GOU process which allows for long memory modeling by the use of a
FBM. 

In order to define a generalized Ornstein-Uhlenbeck process
we define a two-sided
L\'evy process as 
\begin{eqnarray}
\label{eq:two-side-levy}
\xi_t:=\left\{
\begin{array}{ll}
\xi^1_t & \mbox{if}\quad t\ge0  \\
-\xi^2_{-t-} &  \mbox{if}\quad t< 0,
\end{array}
\right.
\end{eqnarray} 
where $\{\xi^1_t\}_{t \ge 0}$ and $\{\xi^2_t\}_{t \ge 0}$ are
independent copies of $\{\xi_t\}_{t\ge 0}$. 
We work throughout with a bivariate complete probability space 
\begin{equation}
\label{eq:omega:space}
\left(\Omega:=\Omega_1\times\Omega_2,\mathscr{F}:=\mathscr{F}_1\otimes\mathscr{F}_2,P:=P_1\otimes
P_2\right). 
\end{equation}
Let $\{\xi_t\}_{t\in \mathbb{R}}$ defined on $\left(\Omega_1,\mathscr{F}_1,
P_1\right) $ be a L\'evy process and a FBM 
$\{B_t^H\}_{t\in\mathbb{R}}$ with Hurst index $H\in(0,1)$ defined on
$\left(\Omega_2 ,\mathscr{F}_2, P_2\right)$ which is
independent of $\{\xi_t\}_{t\in\mathbb{R}}$. A generalized
fractional Ornstein-Uhlenbeck (GFOU) process with
initial value $Y_0 \in L^1(\Omega)$ is defined as 
\begin{eqnarray}
\label{eq:def:GFOU}
Y_t:=e^{-\xi_t}\left(Y_0+\int_0^te^{\xi_{s-}}dB_s^H \right),\quad t\ge0.
\end{eqnarray}
If the initial variable satisfies
\begin{equation}
\label{eq:def:stationary:initial:GFOU}
Y_0 =\int^0_{-\infty} e^{\xi_{s-}} dB_s^H,
\end{equation}
then, for convenience, we sometimes replace $Y=\{Y_t\}_{t\ge0}$ with
\begin{eqnarray}
\label{eq:def:stationary:GFOU}
\overline{Y}_t:= e^{-\xi_t} \int^t_{-\infty} e^{\xi_{s-}} dB_s^H.
\end{eqnarray}
The process $Y$ is regarded as
an extension of $V$ given in (\ref{eq:def:GOU}) where the stochastic process of
integration $\{\eta_t\}_{t\ge0}$ is replaced with a $\{B_t^H\}_{t\ge0}$ with $H\in(0,1)$
and also is regarded as an extended version of the FOU process  
where the integrand is replaced by the exponential of an independent L\'evy process $\xi_t$.
We should remark that $Y$ has jumps caused by the process
$e^{-\xi_t}$.  

The paper is organized as follows. 
In Section \ref{sec:levy} we recall the definition of L\'evy processes and 
summarize properties needed. In Sections
\ref{subsec:Riemann-Stieltjes:p-variation} and \ref{sub:integral:l^2:fbm}
we review Riemann-Stieltjes integrals for functions with bounded $p$-variation and 
the stochastic integral in the $L^2(\Omega)$-sense respectively. In Section
\ref{sec:Existence-Stationality} we investigate the 
existence of the integral in the GFOU process in order to justify the
definition of the GFOU process. 
The stationarity condition and the second order behavior of the GFOU process are discussed in Section
\ref{sec:stationarity:secondorder:GFO}, and we observe the long memory
property. Here we also examine stochastic integrals constructed by a
single FBM, where $\xi$ in the process $Y$ is replaced with $B^H$ used as the
integrator. In Section \ref{sec:sde:gfou} we obtain a stochastic
differential equation, whose solution is given in form of the GFOU process.

We use the following notations throughout. Write $\stackrel{a.s.}{=}$ if
equality holds almost surely. We will take the expectations
for a bivariate process $\{(Z_t^1,Z_t^2)\}_{t\in\mathbb{R}}$. If the expectations only for a process $\{Z^1_t\}_{t\in\mathbb{R}}$ is
considered, we write its expectation as $E_{Z^1}$. 

\section{Preliminaries}
\subsection{L\'evy processes}
\label{sec:levy}
In this subsection we introduce the setup for the L\'evy process.
Let $\xi:=\{\xi_t\}_{t\ge0}$ be a L\'evy process on $\mathbb{R}$ with
$(a_\xi,\nu_\xi,\gamma_\xi)$ generating 
triplet, where $a_\xi\ge0$ and $\gamma_\xi\in\mathbb{R}$ are constants
and a measure $\nu_\xi$ on $\mathbb{R} \setminus\{0\}$ satisfies
$$
\int_{\mathbb{R}\setminus \{0\}}(1\wedge |x|^2)\nu_\xi(dx)<\infty.
$$
We call $\nu_\xi$ the L\'evy measure of
$\xi$.
Then, the characteristic function of $\xi_t$ at time $t=1$ is written as 
\begin{equation}
Ee^{iz\xi_1} =\exp\left[
-\frac{a_\xi}{2} z^2+i\gamma_\xi z+\int_{ \mathbb{R} \setminus
\{0\}}\left(
e^{izx}-1-izx 1_{\{|x|\le 1\}}
\right)\nu_\xi(dx)
\right],\quad z\in\mathbb{R}.
\end{equation}
For more on L\'evy processes and their properties, we refer to
\cite{Sato:1999}. In later sections we consider several examples
related the $\alpha$-stable L\'evy motion with index $0<\alpha<2$, denoted
by $\xi^\alpha:=\{\xi_t^\alpha\}_{t\ge 0}$. It is a L\'evy process and its
generating triplet is $(0,\nu_\alpha,\gamma_\alpha)$ where 
\begin{equation*}
\nu_\alpha(dx):= \left \{
\begin{array}{ll}
c_1x^{-1-\alpha}dx & \text{on}\quad (0,\infty) \\ 
c_2|x|^{-1-\alpha}dx & \text{on}\quad (-\infty,0)
\end{array}
\right.
\end{equation*}
with $c_1\ge0,\ c_2\ge 0$ and $c_1+c_2\ge 0$.

In order to define the 
in (\ref{eq:def:GFOU}), the variation of the L\'evy process plays an
important role. We give a
brief summary based mainly on Section 5.4 of \cite{Dudley:Norvaisa:1998} and
p.408 of \cite{Mikosch:Norvaisa:2000}. Define the $p$-variation for
$0<p<\infty$ of a process $X:=\{X_t\}_{t\in\mathbb{R}}$ on $[t_1,t_2]$
for $t_1<t_2$ in $\mathbb{R}$ as 
\begin{equation}
\label{eq:def-variation}
v_p(X):=v_p(X,[t_1,t_2]):=\sup_{\Delta}\sum_{i=1}^n \left|
X_{s_i}-X_{s_{i-1}}
\right|^p,
\end{equation}
where $\Delta$ is a partition $t_1=s_0<s_1<\cdots<s_n=t_2$ of
$[t_1,t_2]$ and $n\ge1$. If $v_p(X,[t_1,t_2])<\infty$ for all $t_1<t_2$, we say
$X$ has bounded $p$-variation, and if
$v_1(X,[t_1,t_2])<\infty$ we say it is of bounded variation.
 Since every L\'evy process
is a semimartingale $v_p(\xi)<\infty$ for $p\ge2$ (see
\cite{Lepingle:1976}). We will state three useful results which characterize $p$-variation of
L\'evy process in terms of the L\'evy measure. 
Unfortunately we can not find a result which
uniformly characterize the variation in terms of the L\'evy measure.
Assumptions and results are somewhat different from paper to
paper.
The first one is a well-known
result (e.g. Theorem 21,9 (i) of \cite{Sato:1999}).
\begin{quote}
{\bf 1. Bounded variation}

A L\'evy process $\{\xi_t\}_{t\ge0}$ is of bounded variation if and only
if 
$$
a_\xi=0\ \text{and}\ \int_{\mathbb{R}\setminus \{0\}}(1\wedge |x|)\nu_\xi(dx)<\infty.
$$
\end{quote}
The next result is a combination of Theorems 4.1 and 4.2 of \cite{Blumenthal:Getoor:1961} and Theorem 2 of \cite{Monroe:1972}. 
\begin{quote}
{\bf 2.(a) $p$-variation of L\'evy processes} 

Define $\beta=\inf\{\alpha>0:\int_{|x|<1}|x|^\alpha \nu_\xi(dx)<\infty
 \}$. We call $\beta$ the Blumenthal and Getoor index. If the L\'evy
 process
 $\{\xi_t\}_{t\ge0}$ has no Gaussian component $(a_\xi=0)$, then
 \begin{eqnarray}
 \begin{array}{ll}
 v_p(\xi;[0,1])<\infty\quad a.s. & \mbox{if}\quad p>\beta  \\
 v_p(\xi;[0,1])=\infty\quad a.s.  & \mbox{if}\quad p<\beta.
 \end{array}
 \end{eqnarray} 
 \end{quote}
The result by \cite{Bretabnolle:1972} is a sharpened version of {\bf 2.(a)}
but with zero mean assumption.
\begin{quote}
{\bf 2.(b) $p$-variation of L\'evy processes}

Let $1<p<2$ and $\{\xi_t\}_{t\ge0}$ be a L\'evy process
 without Gaussian component $(a_\xi=0)$. Then $v_p(\xi;[0,1])<\infty\ a.s.$ if and only if 
\begin{equation*}
\int_{\mathbb{R}\setminus\{0\}}(1\wedge |x|^p)\nu_\xi(dx)<\infty.
\end{equation*}
Otherwise $v_p(\xi;[0,1])=\infty\ a.s.$
\end{quote}
In particular for 
$\alpha$-stable L\'evy processes we
have the result by \cite{fristedt:taylor:1973} which was stated in
convenient form in \cite{Mikosch:Norvaisa:2000}.  
\begin{quote}
{\bf 3. $p$-variation of $\alpha$-stable L\'evy processes}

Let $\{\xi_t^\alpha\}_{t\ge0}$ be $\alpha$-stable L\'evy motion. Assume
 that $\gamma_\alpha=0$ for $\alpha<1$ and that the L\'evy measure is
 symmetric for $\alpha=1$. Then $v_p(\xi_\alpha)$ is finite or infinite
 with probability $1$ according as $p>\alpha$ or $p\le\alpha$. 
\end{quote}

For the existence of the infinite interval integral in
$\{\overline{Y}_t\}_{t\ge0}$ given in (\ref{eq:def:stationary:GFOU}) we
further need the behavior of $\{\xi_t\}_{t\ge0}$ as $t\to\infty$. Our assumption is that
$\lim_{t\to\infty}\xi_t\stackrel{a.s.}{=}+\infty$.
\cite{doney:maller:2002}[Theorem 4.4] have obtained an equivalent
condition of this in terms of the L\'evy measure
$\nu_\xi$. Since several papers well explain equivalent conditions,
(see p.72 of \cite{erickson:maller:2004} or p.1704 of \cite{lindner:maller:2005})
we do not mention it. Actually if $\lim_{t\to\infty}\xi_t\stackrel{a.s.}{=}+\infty$
holds, we can assert a stronger result, which is more useful for our aim.
\begin{lem}
\label{lem:levy-infty}
Suppose $\lim_{t\to\infty}\xi_t\stackrel{a.s.}{=}+\infty$. Then for
 almost all $\omega_1\in \Omega_1$ there exist
 $\delta>0$ and $t_0=t_0(\omega_1)<\infty$ such that
\begin{equation*}
\xi_t>\delta t\quad \text{for}\quad t\ge t_0.
\end{equation*}
\end{lem}
The proof is a combination of Theorem 4.3 and 4.4 in 
\cite{doney:maller:2002}. Concerning the integral $\int_0^t
e^{-\xi_{s-}}d\eta_s$ in the GFOU process given by (\ref{eq:def:GOU}), \cite{erickson:maller:2004}
have characterized the convergence of improper integral $\int_0^\infty
e^{-\xi_{s-}}d\eta_s$ in terms of the L\'evy measure of
$\{(\xi_t,\eta_t)\}_{t\ge0}$, in which the condition 
$\lim_{t\to\infty}\xi_t\stackrel{a.s.}{=}+\infty$ was used.

\subsection{Integrals with respect to functions with unbounded variation} 
\label{sec:integrals}
\subsubsection{Riemann-Stieltjes integrals with $p$-variation}
\label{subsec:Riemann-Stieltjes:p-variation}
We review several useful definitions of integrals of functions which have
unbounded variation but bounded $p$-variation. The excellent
introduction to this area is given by \cite{Dudley:Norvaisa:1998}. Let
$f$ and $g$ be real functions on $[a,b]$. Define
$\kappa=\{u_1,\ldots,u_n\}$ to be an intermediate partition of $\Delta=[a=s_0<s_1<\cdots<s_n=b]$
given as in (\ref{eq:def-variation}), namely, $s_{j-1}\le u_i \le s_j$ for
$i=1,\ldots,n$. A Riemann-Stieltjes sum is defined as 
$$
S_{RS}(f,g,\Delta,\kappa):=\sum_{i=1}^{n}f(u_i)\left[
g(s_i)-g(s_{i-1})
\right].
$$
Then we say that the Riemann Stieltjes integral 
exists and
equals to $I$, if for every $\epsilon>0$, there exists $\delta>0$ such
that
$$
\left| S_{RS}(f,g,\Delta,\kappa)-I
\right|<\epsilon
$$
for all partition $\Delta$ with mesh $\max\ (s_i-s_{i-1})<\delta$ and for
all intermediate partitions $\kappa$ of $\Delta$. The following theorem
was proved by \cite{Young:1936}. (See also Theorem 2.4 of
\cite{Mikosch:Norvaisa:2000} or Theorem 4.26 of
\cite{Dudley:Norvaisa:1998}.) 
\begin{thm}
\label{thm:rmintegra:p-variation}
Assume $f$ has bounded $p$-variation and $g$ has bounded $q$-variation
 on $[a,b]$ for some $p,q>0$ with $p^{-1}+q^{-1}>1$. Then the integral
 $\int_a^bf dg$ exists in the Riemann-Stieltjes sense, and the inequality 
$$
\left|\int_a^b f dg \right|\le 
C_{p,q}
(v_p(f))^{1/p}(v_q(g))^{1/q}
$$
holds with $C_{p,q}=\sum_{n\ge1}n^{-(p^{-1}+q^{-1})}$. 
\end{thm}

\subsubsection{Integral with respect to FBM in $L^2$-sense}
\label{sub:integral:l^2:fbm}
Another definition of the integral is in $L^2(\Omega_2)$-sense. 
Stochastic integrals with respect to FBM is sometimes defined as the
$L^2(\Omega_2)$-limits of integrals of step functions (see e.g. \cite{lin:1995}).
We see this when
$B^H$ is the Brownian motion $B^{1/2}$. If a function
$f(u):\mathbb{R}\to\mathbb{R}$ satisfies $f(u)\in L^2(\mathbb{R})$, there exist step functions
$$
f_n(u):=\sum_{i=1}^n f_i 1_{\{u_i<u\le u_{i+1}\}},\quad
-\infty<u_1<\ldots<u_{n+1}<\infty, \quad f_i\in \mathbb{R},\ n\in \mathbb{N}. 
$$
such that $f_n$ converges to $f$ in $L^2(\mathbb{R})$. Then
the integral $\int_{\mathbb{R}}f(u)dB_u^{1/2}$ is the 
$L^2(\Omega_2)$-limit of the integrals of step functions
$$
\int_{\mathbb{R}}f_n(u)dB_u^{1/2}=\sum_{i=1}^n f_i \left(
B_{u_{i+1}}^{1/2}-B_{u_i}^{1/2}
\right),
$$
since $L^2(\mathbb{R})$ and $L^2(\Omega)$ are isometric and their inner
products are equal, namely,
$$
E\left[
\left(\int_{\mathbb{R}}
(f_n(u)-f(u))dB_u^{1/2}
\right)^2
\right]=\int_\mathbb{R}\left(f_n(u)-f(u)\right)^2du.
$$

\cite{pipiras:taqq:2000b} have investigated a similar characterization
for the integral of $B^H$ when $H\neq \frac{1}{2}$. We apply this to the 
existence of the improper integral in the GFOU process $\{\overline{Y}_t\}$ afterward. Define the linear space
$$
\overline{Sp}(B^H):=\left\{ X:\int_\mathbb{R}f_n(u)dB^H_u
\to
X\ \text{in}\ L^2(\Omega) \ \text{for some sequence}\ (f_n)_{n\in \mathbb{N}}\ \text{(step functions)} \right\}.
$$
\cite{pipiras:taqq:2000b} have analyzed the functional space of the
integrand $f(u)$ in which it can be asymptotically approximated
by $f_n(u)$ and $\int_{\mathbb{R}}f(u)dB_u^H$ is well-defined. 
For $H\in(0,\frac{1}{2})$ they succeeded in specified a Hilbert space of
functions on the real line which is isometric to $\overline{Sp}(B^H)$.
However, for $H\in (\frac{1}{2},1)$ they had difficulty in
finding the corresponding isometric space, and as second best
they analyzed inner product spaces in which the integral
with $B^H$ ($H\in (\frac{1}{2},1)$) is well-defined. We give only one such
inner product space and its inner product for $B^H$ with $H\in(\frac{1}{2},1)$.
Other inner product spaces do not seem to work for our purpose since they
require e.g. characteristic function of $f$ or fractional derivative of
$f$ which do not exists in our case where $f(u)=e^{\xi_u}$. (See Section 7 of \cite{pipiras:taqq:2000b}.) 
The space is 
\begin{equation}
\label{eq:innerproduct:space}
\left|\Lambda\right|^H = \left\{ 
f:
 \langle |f|,|f|\rangle_{|\Lambda|^H} <\infty\right\},\quad 
H\in\left(\frac{1}{2},1\right),
\end{equation}
where
\begin{equation*}
\langle f,g \rangle_{|\Lambda|^H}=c_H \int_{\mathbb{R}}\int_{\mathbb{R}}f(u)g(v)|u-v|^{2H-2}dudv
\end{equation*}
is the inner product with $c_H:=H(2H-1)$.

\section{Existence of the integral}
\label{sec:Existence-Stationality}
In this Section we analyze the existence of the integral in the GFOU
process given in (\ref{eq:def:GFOU}).
The definition of the GFOU process includes an integral with respect to
FBM. Since paths of FBM are of infinite variation and since FBM with
$H\neq\frac{1}{2}$ is not a semimartingale, the stochastic integral with
respect to FBM $(H\neq \frac{1}{2})$ is not an It\^o integral. 
Additionally, the integrand of the GFOU process is random and the infinite interval
integral (\ref{eq:def:stationary:initial:GFOU}) is needed for its stationarity. 
We apply two approaches of the integral in Section \ref{sec:integrals} in order to cope with
these problems.  

\begin{prop}
\label{prop:existence:pathwisely}
Let $B^H:=\{B_t^H\}_{t\in\mathbb{R}}$ be a FBM with $H\in(0,1)$ and
 $\xi:=\{\xi_t\}_{t\in\mathbb{R}}$ be an independent two-sided L\'evy
 process. Assume that $\{\xi_t\}_{t\in\mathbb{R}}$ has bounded
 $p$-variation for $0<p<\infty$. Then the integral
 $\int_0^te^{\xi_{s-}}dB_s^H$, $0<t<\infty$ exists in the path-wise Riemann-Stieltjes
 sense if $\frac{1}{p}+H>1$.
 Furthermore with $q>\frac{1}{H}$ and $C_{p,q}=\sum_{n\ge
 1}n^{-(p^{-1}+q^{-1})}$ we have 
 \begin{equation*}
 \left|
\int_0^t e^{\xi_{s-}}dB_s^H
 \right| \le C_{p,q} \bigg( \sup_{s\in [0,t]}e^{\xi_s} \bigg)
 \left(v_p\left(\xi \right)\right)^{1/p}\left(v_q\left(B^H \right)\right)^{1/q}\quad P-a.s.
 \end{equation*}
\end{prop}
The proof of Proposition \ref{prop:existence:pathwisely} is obvious from
the continuity of the exponential function and we omit it.
\begin{rem}
If $p=1$ in Proposition \ref{prop:existence:pathwisely}, $\xi$ has bounded variation and we can
 define path-wise integrals for all $B^H$ with $H\in(0,1)$. On the other hand
 if $H>\frac{1}{2}$, we can define the path-wise integral for any integrands
 $\xi$ since $v_p(\xi)<\infty$ for $p>2$.  
\end{rem}

The reason why we need the path-wise definition besides the
$L^2(\Omega)$p-approach is that we want to define the integral for
$H\in(0,\frac{1}{2})$, which we could not obtain in the
$L^2(\Omega)$-approach. It also gives several useful tools easily, such as integration by parts or
chain rule for analyzing stochastic differential equations in Section
\ref{sec:sde:gfou}. \\


\begin{ex}
Set $\xi:=\xi^\alpha$ be an $\alpha$-stable L\'evy motion with index
 $\alpha$ defined in Subsection \ref{sec:levy}. 
 Then under assumption in {\bf 3. $p$-variation of $\alpha$-stable
 L\'evy processes}
 and the assumption $\frac{1}{\alpha}+H>1$. the integral
 $e^{-\xi^\alpha_t}\int_0^te^{\xi^\alpha_{u-}}dB_u^H $ 
 exists in the path-wise sense. Hence the GFOU process is well-defined. 
 Both $\xi^\alpha$
 and $B^H$ have infinite variation
 and the former is known to be extremely heavy tailed process. 
\end{ex}

Recall that we write $Y_t=Y_t(\omega_1,\omega_2)$ to
emphasis that it is a function from probability space
(\ref{eq:omega:space}) and that we write $\overline{Y}_t$ if initial
value $Y_0$ satisfies (\ref{eq:def:stationary:initial:GFOU}).
Define approximating step functions of the integrand
$e^{-\xi_t+\xi_{u-}}$ of $\overline{Y}_t$ as 
\begin{equation*}
f_{t,n}(u;\omega_1):= \sum_{i=1}^n f_i(u;\omega_1) 1_{\{u_{i-1}^n<u\le u_i^n\}},
\end{equation*}
where 
$$
f_i(u;\omega_1):=e^{-\xi_t(\omega_1)+\xi_{u_{i-1}^n}(\omega_1)}. 
$$
%
Here $u_i^n$'s are points in $[-N_n,t]$ such that $-N_n = u_0^n<u_1^n<\cdots<u_n^n=t$ and as $n\to\infty$
$\max\ (u_i^n-u_{i-1}^n)\downarrow 0$ and $N_n\uparrow \infty$. By
integrating $f_{t,n}$ with respect to $B_t^H$, we also define
approximating sequence as 
\begin{eqnarray*}
Z_t^n(\omega_1,\omega_2) &:=& \int_{-\infty}^t
 f_{t,n}(u;\omega_1)dB^H_u(\omega_2) \nonumber \\
&=&\label{eq:func:approxi:step}
\sum_{i=1}^n f_{i}(u;\omega_1)(B_{u_i^n(\omega_2)}^H-B_{u_{i-1}^n}^H(\omega_2)).
\end{eqnarray*}
In the following theorem we define the integral $\int_{-\infty}^t
e^{-\xi_t-\xi_{u-}}dB_u^H$ as the limit in probability of the 
$Z_t^n$ as $n\to \infty$. 
The reason why we need this approach is that with only 
path-wise definitions we find difficulty to treat improper
integrals. For the existence of the improper integral we should consider
long time 
($t\to\infty$) behavior of both $\{B_t^H\}_{t\in\mathbb{R}}$ and
$\{\xi_t\}_{t\in \mathbb{R}}$ path-wisely, which seems to be not an easy task. 
Additionally, this approach is well-matched with the analysis of the
second order behavior.

\begin{thm}
\label{thm:GFOU:existence:fubini}
Let $\{B_t^H\}_{t\in\mathbb{R}}$ be a FBM with $H\in(\frac{1}{2},1)$ and
 $\{\xi_t\}_{t\in\mathbb{R}}$ be an independent two-sided L\'evy
 process. If $\lim_{t\to\infty}\xi_t \stackrel{a.s.}{=} +\infty$, 
 then for each $t\ge0$ $Z_t^n$ given in (\ref{eq:func:approxi:step})
 converges in probability to a limit defined by $\overline{Y}_t$ and
 which does not depend on the sequence $u_i^n$. 
 If  further $E[e^{-2\xi_1}]<1$, then $Z_t^n$ converges to $\overline{Y}_t$ in $L^2(\Omega)$
 and it follows that 
 \begin{equation*}
 E[\overline{Y}_t]=E\left[
 \int_{-\infty}^t e^{-(\xi_t-\xi_{u-})}dB_u^H
 \right]=0
 \end{equation*}
 and for $0<s\le t$,  
 \begin{eqnarray}
 E[\overline{Y}_s \overline{Y}_t] 
 &=& 
 \int_{-\infty}^s \int_{-\infty}^t  
 E_\xi[e^{-\xi_s+\xi_{u-}-\xi_t+\xi_{v-}}]
 |u-v|^{2H-2}dudv. \label{thm:GFOU:eq:second:}
 \end{eqnarray}
\end{thm}

\begin{rem}
\label{rem:integral2}
 (a) For $H\in (0,\frac{1}{2})$ we could not validate the existence of
 the improper integral in $\overline{Y}_t$.\\
 (b) Under the only condition that $\lim_{t\to\infty}\xi_t \stackrel{a.s.}{=}+\infty$,
 we cannot prove the $L^2(\Omega)$ convergence. \\
 (c) The condition $E[e^{-2\xi_1}]<1$ implies
  $\lim_{t\to\infty}\xi_t \stackrel{a.s.}{=} +\infty$, which is shown in
  Proposition 4.1 of \cite{lindner:maller:2005}.
\end{rem}

\noindent
{\it Proof of Theorem \ref{thm:GFOU:existence:fubini}}. \\
We check that for each $\omega_1$ the integrand
$e^{-\xi_t(\omega_1)+\xi_{u-}(\omega_1)}\in |\Lambda|^H$ given in
(\ref{eq:innerproduct:space}). Since $\xi_t(\omega_1)$ is constant, we
drop it and only show $e^{\xi_{u-}(\omega_1)}\in |\Lambda|^H$. 
The function
$e^{\xi_{u-}(\omega_1)}$ is c\`agl\`ad and bounded on
any finite interval. Due to Lemma \ref{lem:levy-infty} and to the symmetry of
two-sided L\'evy processes there exists $T(\omega_1)<0$ such that for all $u\le T(\omega_1)$, 
$\xi_{u-}<\delta u$ where $\delta$ is some positive constant. Then for
each $\omega_1$, 
\begin{eqnarray*}
&& \int_{-\infty}^t\int_{-\infty}^t
 e^{\xi_{u-}(\omega_1)+\xi_{v-}(\omega_1)}|u-v|^{2H-2}dudv <\infty
\end{eqnarray*}
 is obvious.
 Hence we can utilize $L^2(\Omega_2)$-integral theory in Section
 \ref{sub:integral:l^2:fbm}. Namely, for each $t\ge0$ and for each
 fixed $\omega_1$, $Z_t^n(\omega_1,\cdot)$ converges
 in $L^2(\Omega_2,P_2)$. 
 Moreover $Z_t^n$ converges in probability on $(\Omega,P)$ for each $t\ge0$
 since sequence $Z_t^n$ satisfies the Cauchy criterion, as seen by
 \begin{eqnarray*}
 && \lim_{n,m \to\infty}P\left(\left|Z_t^n-Z_t^m \right|>\epsilon
		      \right) \\
 && \quad = \lim_{n,m \to\infty}\int_{\Omega_1}\int_{\Omega_2} 1_{\{|Z_t^n-Z_t^m|>\epsilon
 \}}dP_2dP_1 \\
 && \quad =
 \int_{\Omega_1}\lim_{n,m
 \to\infty}P_2\left(\left|Z_t^n(\omega_1)-Z_t^m(\omega_1)\right|>\epsilon \right)dP_1 =0.
 \end{eqnarray*}
 The limit is called $\overline{Y}_t$ and it is $\mathscr{F}_1\otimes\mathscr{F}_2$
 measurable for each $t$.

 Now, with $E[e^{-2\xi_1}]<1$, we prove the $L^2(\Omega)$-convergence. 
 We have $E[e^{-\xi_1}]<1$ as well, hence    
$$
\theta_1:=-\log E[e^{-\xi_1}]>0\quad \text{and}\quad \theta_2:=-\log E[e^{-2\xi_1}]>0.
$$
Then using the covariance of the FBM, we have
\begin{eqnarray}
E[(Z_t^n)^2] &=& E\left[
\sum_{i=1}^n\sum_{j=1}^ne^{-2\xi_t+\xi_{u_{i-1}^n}+\xi_{u_{j-1}^n}}(B^H_{u_i^n}-B^H_{u_{i-1}^n})(B^H_{u_j^n}-B^H_{u_{j-1}^n})
\right] \noindent \nonumber \\
&=& \sum_{i=1}^n \sum_{j=1}^n \bigg\{
 E_\xi \left[ e^{-2\xi_t+\xi_{u_{i-1}^n}+\xi_{u_{j-1}^n}} \right]\frac{1}{2}\big(
 -|u_i^n-u_j^n|^{2H}+|u_{i-1}^n-u_j^n|^{2H} \nonumber \\
&& \hspace{2cm} +|u_{i}^n-u_{j-1}^n|^{2H}-|u_{i-1}^n-u_{j-1}^n|^{2H}
 \big)
 \bigg\} \nonumber \\
&=& c_H\int_{-\infty}^t\int_{-\infty}^t \Gamma^n(u,v)|u-v|^{2H-2}dvdw,
\label{eq:int:approx}
\end{eqnarray}
where 
\begin{eqnarray*}
\Gamma^n_t(u,v) &=& \sum_{i=1}^n\sum_{j=1}^n
\bigg\{
1_{\{u_{i-1}^n<u\le u_i^n, u_{j-1}^n<v\le
u_j^n, i\ge j\}}
e^{-\theta_2(t-u_{i-1}^n)-\theta_1(u_{i-1}^n-u_{j-1}^n)} \\
&& \hspace{1.3cm}
+1_{\{u_{i-1}^n<u\le u_i^n, u_{j-1}^n<v\le
u_j^n, i< j\}}
e^{-\theta_2(t-u_{j-1}^n)-\theta_1(u_{j-1}^n-u_{i-1}^n)}
\bigg\},
\end{eqnarray*}
which obviously converges point-wise to 
\begin{eqnarray*}
 1_{\{u\ge
 v\}}e^{-\theta_2(t-u)-\theta_1(u-v)} + 1_{\{u<
 v \}}e^{-\theta_2(t-v)-\theta_1(v-u)}.
\end{eqnarray*}
Since $|\Gamma^n_t(v,w)|<M'$ for some $M'>0$ uniformly in $n,\ t\le 0$, 
(\ref{eq:int:approx}) is bounded by $c_H M' t^H$. Furthermore according
to usual Fubini's theorem, it also
follows that 
\begin{eqnarray*}
E\left[(\overline{Y}_t)^2\right] 
&=& \int_{\Omega_1}\int_{\Omega_2}\left(
\overline{Y}_t(\omega_1,\omega_2) 
\right)^2 dP_2dP_1 \\
&=& c_H E_\xi \left[
\int_{-\infty}^t \int_{-\infty}^t e^{-2\xi_t+\xi_{u-}+\xi_{v-}}
|u-v|^{2H-2}dudv 
\right] \\
&=& 2c_H\int_{-\infty}^t\int_{-\infty}^t e^{-\theta_2(t-u)-\theta_1(v-u)}
1_{\{u\ge v\}}|u-v|^{2H-2}dudv<\infty. 
\end{eqnarray*}
Observe that this integral is finite.
Accordingly
$E[(Z_t^n)^2]\to E\left[(\overline{Y}_t)^2\right]$ as $n\to\infty$.
Now we can apply Theorem 4.5.4 of \cite{Chung:2001} to
$Z_t^n$ and obtain the $L^2(\Omega)$-convergence. In consequence
$E[\overline{Y}_s\overline{Y}_t]$ turns out to be finite and
equation (\ref{thm:GFOU:eq:second:}) 
follows from Fubini's theorem. 
Finally $E[Z_t^n]=0,\ n\in\mathbb{N}$
implies $E[\overline{Y}_t]=0.$
Hence the proof is
complete. \hfill $\Box$\\

The process $\{\overline{Y}_t\}_{t\in\mathbb{R}}$ obtained in Theorem
\ref{thm:GFOU:existence:fubini} is the GFOU process with initial value $\overline{Y}_0$.
 We close this section with the following concluding Remarks.
\begin{rem}
\label{rem:convergence:integral}
 (a)
 In both Proposition \ref{prop:existence:pathwisely} and Theorem
 \ref{thm:GFOU:existence:fubini} $\xi$
 is independent of $B^H$ 
 and we have 
 $$
 e^{-\xi_t}\int_{a}^{b} e^{\xi_{u-}}dB_u^H =
 \int_{a}^{b}e^{\xi_t+\xi_{u-}}dB_u^H, \quad -\infty \le a<b<\infty.
 $$
 This is not allowed in usual theory of stochastic integrals related to semimartingale
 (\cite{Protter:2004}) since $\{\xi_t\}_{t\in\mathbb{R}}$ is not adapted. \\
 (b) From Theorem \ref{thm:GFOU:existence:fubini} $\lim_{t\to\infty}\xi_t \stackrel{a.s.}{=}+\infty$
  implies $|Y_0|=\left|\int^0_{-\infty} e^{\xi_{s-}} dB_s^H \right|<\infty\ a.s.\ \Omega.$ Therefore together with
  Proposition \ref{prop:existence:pathwisely} we can also treat
  $\{\overline{Y}_t\}_{t\ge0}=\{Y_t+Y_0\}_{t\ge0}$ path-wisely. \\
 (c) We should mention that \cite{erickson:maller:2004}[p.81] gave  
  another idea for improper integrals of a function of L\'evy processes
  with respect to FBM:$\int_0^\infty g(\xi_t)dB_t^H$. Investigation of their idea is also our next concern.
\end{rem}

\section{Stationarity and Second order behavior of GFOU processes}
\label{sec:stationarity:secondorder:GFO}
Here we investigate the strict stationarity and the second order
behavior of the GFOU process $\overline{Y}:=\{\overline{Y}_t\}_{t\ge0}$. Since we could not
validate the existence of $\overline{Y}$ with Hurst parameter
$H\in(0,\frac{1}{2})$ and 
since our main concern in this paper is the long memory
case, we confine our results to the case $H\in(\frac{1}{2},1)$ throughout
this section.

The stationarity of $\overline{Y}$ is as follows.
\begin{prop}
\label{prop:stationarity-GFOU}
 If $\lim_{t\to\infty}\xi_t\stackrel{a.s.}{=}+\infty$, then
 $\overline{Y}_t$ exists for all $t$ and the process $\overline{Y}:=\{\overline{Y}_t\}_{t\ge 0}$ is strictly stationary.
\end{prop}

\noindent
{\it Proof of Proposition \ref{prop:stationarity-GFOU}}. \\
Let $0\le t_1<t_2<\ldots<t_m<\infty,\ m\in \mathbb{N}$ and $h>0$. We use
the sequence $Z_t^n$ given in (\ref{eq:func:approxi:step}). Since both
$\{\xi_t\}_{t \in \mathbb{R}}$ and $\{B_t^H\}_{t\in\mathbb{R}}$ have stationary increments, does the pair
$\{(\xi_t,B_t^H)\}_{t\in\mathbb{R}}$ as well because of independence. Thus,
\begin{eqnarray*}
Z_{t_i}^n &=& \sum_{j=1}^n
 e^{-\xi_t+\xi_{u_{j-1}^n}}(B_{u_j^n}^H-B_{u_{j-1}^n}^H),\quad -N_n=u_1^n<\ldots<u_n^n=t_i  \\
&\stackrel{d}{=}& \sum_{j=1}^n
 e^{-\xi_{t+h}+\xi_{u_{j-1}^n+h}}(B_{u_j^n+h}^H-B_{u_{j-1}^n+h}^H) \\
&=& Z_{t_{i-1}+h}^n,\quad 1\le i\le m,
\end{eqnarray*}
simultaneously in $i$ and we have 
$$
(Z_{t_1}^n,\ldots,Z_{t_m}^n) \stackrel{d}{=} (Z_{t_1+h}^n,\ldots,Z_{t_m+h}^n).
$$
Then by virtue of Theorem \ref{thm:GFOU:existence:fubini} as $n\to\infty$, $(Z_{t_1}^n,\ldots,Z_{t_m}^n)$
converges in probability to $(Y_{t_1},\ldots,Y_{t_n})$ and 
$(Z_{t_1+h}^n,\ldots,Z_{t_m+h}^n)$ converges in probability to
$(Y_{t_1+h},\ldots,Y_{t_n+h})$. 
This yields
\begin{equation*}
(Y_{t_1},\ldots,Y_{t_n})\stackrel{d}{=} (Y_{t_1+h},\ldots,Y_{t_n+h})
\end{equation*}
and the conclusion holds. \hfill $\Box$

\begin{rem}
(a) The logic in the proof works even in the case
 $H\in(0,\frac{1}{2})$. If the assumptions of Proposition
 \ref{prop:existence:pathwisely} are satisfied and if
 the integral $Y_0=\int_{-\infty}^0e^{\xi_{s-}}dB_s^H$ with
 $H\in(0,\frac{1}{2})$ exists, then
 $\overline{Y}$ with $H\in(0,\frac{1}{2})$ is defined
 and is strictly stationary. \\
 (b) In connection with the GOU process $\{V\}_{t\ge0}$ given in
 $(\ref{eq:def:GOU})$. Proposition \ref{prop:stationarity-GFOU}
 corresponds to Theorem 3.1 of \cite{carmona:petit;yor:2001}, where
 $\xi$ and $\eta$ are independent, and under conditions of $a.s.$
 convergence of the integral $\int_0^\infty e^{-\xi_{s-}}d\eta_s$ and 
 $\lim_{t\to \infty}\xi_t\stackrel{a.s.}{=}+\infty$ the stationary
 version exists and equals in distribution to $V_\infty$. 
\end{rem}

Next we investigate the second order behavior of
$\overline{Y}$ and derive the auto-covariance function
explicitly. What should be remarked is that while the auto-covariance
function of the GOU process $V$ given in $(\ref{eq:def:GOU})$ decreases exponentially (Theorem 4.2
of \cite{lindner:maller:2005}), that of the FOU process $\{X_t^H\}_{t\ge 0}$ given
in $(\ref{eq:def:GOU})$ decays
like a power function (Theorem 2.4 and Corollary 2.5 in
\cite{cheridito:kawaguchi:maejima:2003}). Since
$\overline{Y}$ is regarded as a version of GOU processes
and FOU processes, 
this investigation is interesting. We utilize results in Theorem 3.1 and
obtain Theorem \ref{thm:second-order-behavior-stationary:GFOU} and
Corollary \ref{cor:second-order-behavior:GFOU} below. Note that even the 
existence of $\overline{Y}$ and the equation (\ref{thm:GFOU:eq:second:}) are
obtained several difficulties still lay in calculating the
auto-covariance function. The integrand in (\ref{thm:GFOU:eq:second:})
is regarded as exponential moment of of $4$ dependent random variable,
i.e. $E_\xi[e^{-(\xi_s-\xi_{u-}+\xi_t-\xi_{v-})}]$ and dependencies of
these variables are different in the order of $s,u,t$ and $v$. We also
require that after $E_\xi[e^{-(\xi_s-\xi_{u-}+\xi_t-\xi_{v-})}]$ is
calculated the double integral in (\ref{thm:GFOU:eq:second:}) has a
suitable representation for our purpose. In the proofs of Theorem \ref{thm:second-order-behavior-stationary:GFOU} and
Corollary \ref{cor:second-order-behavior:GFOU} we will see how to get
over these difficulties.  
 Proofs of Theorem
 \ref{thm:second-order-behavior-stationary:GFOU} and Corollary
 \ref{cor:second-order-behavior:GFOU} are given in Appendix
 \ref{sec:pf:secondorder} since 
 they require a lot of technical and tedious calculations. 
\begin{thm}
 \label{thm:second-order-behavior-stationary:GFOU}
Let $\{B_t^H\}_{t\in\mathbb{R}}$ be a FBM with $H\in(\frac{1}{2},1)$ and
 $\{\xi_t\}_{t\in\mathbb{R}}$ be an independent two-sided L\'evy
 process. Suppose that $E[e^{-2\xi_1}]<1$. Then the 
 stationary version $\overline{Y}:=\{\overline{Y}_t\}_{t\ge 0}$ exists and for any
 $s>0,t\ge 0$, we have
 \begin{eqnarray*}
 && \mathrm{Cov}(\overline{Y}_t,\overline{Y}_{t+s})= \mathrm{Cov}(\overline{Y}_0,\overline{Y}_{s})\\
 &&\quad = c_H\bigg(
 \frac{2e^{-\theta_1s}}{\theta_2\theta_1^{2H-1}}\Gamma(2H-1)
 -\frac{e^{-\theta_1s}}{2\theta_1^{2H}}\Gamma(2H-1) \\
 &&\quad\quad + \frac{e^{-\theta_1s}s^{2H-1}}{2\theta_1(2H-1)} {}_1\mathrm{F}_1(2H-1,2H;\theta_1s)
 +\frac{e^{2\theta_1s}}{2\theta_1^{2H}}\Gamma(2H-1,\theta_1s)  \bigg) \\
 &&\quad = c_H \bigg[
 \frac{2e^{-\theta_1s}}{\theta_2\theta_1^{2H-1}}\Gamma(2H-1)
 -\frac{e^{-\theta_1s}}{2\theta_1^{2H}}\Gamma(2H-1) \\
 &&\quad \quad+ \frac{s^{2H-1}}{2\theta_1(2H-1)}\sum_{n=0}^\infty
  \prod_{k=1}^{n+1}(2H-k)\bigg\{
 (\theta_1s)^{-(n+1)}-(-\theta_1 s)^{-(n+1)}\frac{\gamma(n+1,\theta_1s)}{n!}
 \bigg\}
 \bigg],
 \end{eqnarray*}
 where $\theta_1:=-\log (E [ e^{-\xi_1} ] )>0$ and
 $\theta_2:=-\log(E[e^{-2\xi_1}])>0$. Here $\gamma(\cdot,\cdot)$ and
 $\Gamma(\cdot,\cdot)$ are incomplete gamma functions in 8.350 of
 \cite{Gradshteyn:Ryzhik:2000} and ${}_1\mathrm{F}_1(\cdot,\cdot;\cdot)$
 is the confluent hyper-geometric function in 9.210 of
 \cite{Gradshteyn:Ryzhik:2000}. 
\end{thm}
Note that since $\gamma(n+1,s\theta_1)\to \Gamma(n+1)=n!$ as $s\to\infty$,
 we have 
 \begin{eqnarray*}
 \mathrm{Cov}(\overline{Y}_t,\overline{Y}_{t+s}) &=&
  H \sum_{n=1}^\infty \prod_{k=1}^{2n-1} (2H-k) \theta_1^{-2n}s^{2H-2n}+
  O(e^{-\theta_1s})\\
 &=& O(s^{2H-2}).
 \end{eqnarray*}
This conclude that $\overline{Y}$ with $H\in(\frac{1}{2},1)$ is a long memory
process. While we obtained
$\mathrm{Cov}(\overline{Y}_t,\overline{Y}_{t+s})$ using special
functions in Theorem \ref{thm:second-order-behavior-stationary:GFOU}, it
is mainly for numerical purpose since for such functions useful
softwares are available.

Next we investigate long time dependence of $\{Y_t\}_{t\ge 0}$ with
the initial value $Y_0:=X\in L^2(\Omega)$ where $X$ is independent of
$\{\xi_t\}_{t\ge0}$ and $\{B_t^H\}_{t\ge0}$.
\begin{cor}
\label{cor:second-order-behavior:GFOU}
Let $Y:=\{Y_t\}_{t\ge0},\ H\in(\frac{1}{2},1)$ be a GFOU process with the
 initial value $X\in L^2(\Omega)$, where $X$ is independent of
 $\xi:=\{\xi_t\}_{t\ge 0}$ and $B^H:=\{B_t^H\}_{t\ge 0}$. Then for fixed $t\ge 0$ as $s\to\infty$.
 \begin{eqnarray*}
 \mathrm{Cov}(Y_t,Y_{t+s}) &=& H\sum_{n=1}^\infty
  \prod_{k=1}^{2n-1}(2H-k)\theta_1^{-2n}\{
 s^{2H-2n}-e^{-\theta_1t}(t+s)^{2H-2n}
 \}+O(s^{2H-2N-2}) \\
 &=& O(s^{2H-2}).
 \end{eqnarray*}
\end{cor}

We see what happens to the second order behavior of the process
$Y$ if $\xi$ in
$Y$ is replaced with
$B^H$ which is the same process as the variable
of integration. Although we expect long memory this does not
hold.
Note that we need only the probability space
$\left(\Omega_2,\mathscr{F}_2,P_2\right)$ here. For
$H\in(\frac{1}{2},1)$ and the initial random variable $X\in L^2(\Omega_2)$ independent of
$\{B_t^H\}_{t\in\mathbb{R}}$, define
$$
W_t=e^{-B_t^H}\left(X+\int_0^t e^{B_{u-}^H}dB_u^H \right).
$$
To analyze $\{W_t\}_{t\ge0}$ we
use the path-wise integral theory (see Subsection
\ref{subsec:Riemann-Stieltjes:p-variation}).
Let $f$ be continuous differentiable and $F(x)=F(0)+\int_0^x
f(y)dy$. Then with $H\in(\frac{1}{2},1)$ it follows that 
$$
F(B_t^H)-F(B_0^H)=\int_0^t f(B_t^H)dB_u^H\quad a.s.
$$
By setting $f=e^t$ in above we obtain
$$
W_t=1+e^{-B_t^H}(X-1).
$$
\begin{prop}
\label{prop:cov-corr-FBM-FBM}
Let $H\in(\frac{1}{2},1)$ and $t,s\ge0$.
Define $M_1:=\left(E[X]-1\right)^2$ and $M_2:=E[(X-1)^2]$. 
 The process $\{W_t\}_{t\ge0}$ has
 the following auto-covariance and correlation functions. 
\begin{eqnarray*}
\mathrm{Cov}\left(W_t,W_{t+s} \right) &=&
e^{\frac{1}{2}\{
t^{2H}+(t+s)^{2H} \}}
\left\{
M_2\, e^{\frac{1}{2}\{t^{2H}+(t+s)^{2H}-s^{2H}\}}
-M_1
\right\} \\
&=& O\left(e^{\frac{1}{2}\{(t+s)^{2H}+s^{2H-1}\} } \right)\ \text{as}\
 s\to\infty. \\
\mathrm{Corr}\left(W_t,W_{t+s}\right) &=& \frac{
M_2\, e^{\frac{1}{2}(t^{2H}-s^{2H})}- M_1\, 
e^{-\frac{1}{2}(t+s)^{2H}} }
{\sqrt{M_2\, e^{t^H}-M_1}\sqrt{M_2-M_1\, e^{-(t+s)^{2H}}}}\\
&= &O(e^{-\frac{1}{2}s^H})\ \text{as}\ s\to \infty. 
\end{eqnarray*}
\end{prop}
We also consider the drift added process
$$
\widehat{W}_t=e^{-(B_t^H+at)}\left(
X+\int_0^t e^{B_u^H+au}d(B_u^H+au)
\right),
$$
where $a>0$ and $X\in L^2(\Omega_2)$ is independent of
$\{B_t^H\}_{t\ge0}$. Even if a drift is added, usual path-wise integral
works and 
$$
\widehat{W_t}=1+(X-1)e^{-(B_t^H+at)}
$$
holds. Then the auto-covariance and correlation functions of $\{\widehat{W}_t\}_{t\ge0}$
are calculated is a similar manner and become
\begin{eqnarray*}
\mathrm{Cov}(\widehat{W}_t,\widehat{W}_{t+s}) &=&
 e^{-a(2t+s)}\mathrm{Cov}\left(W_t,W_{t+s} \right), \\
\mathrm{Corr}(\widehat{W}_t,\widehat{W}_{t+s}) &=&
\mathrm{Corr}\left(W_t,W_{t+s}\right).
\end{eqnarray*}
Thus our conclusion here is that even if a drift is added it does not
have long memory. 

\section{Stochastic differential equation related with GFOU processes}
\label{sec:sde:gfou}
We analyze a stochastic differential equation of which a solution is given
by the GFOU process. Let $U:=\{U_t\}_{t\ge0}$ be a L\'evy process with generating
triplet $(a_U,\nu_U,\gamma_U)$. Assume that the L\'evy measure
$\nu_U$ has no mass on $(-\infty,-1]$. The Dol\'eans-Dade exponential
of $U_t$ is written as $\mathcal{E}(U_t)=e^{-\xi_t}$ where
\begin{equation*}
\xi_t=-U_t+\frac{a_\xi}{2}t-\sum_{0<s\le t}\left(\log(1+\Delta U_s)-\Delta
				     U_s\right). 
\end{equation*} 
See Section 2.2 of \cite{erickson:maller:2004}. Here $\xi_t$ is the L\'evy
processes. 

\begin{prop}
\label{prop:SDE:GFOU}
Under the assumption in Proposition \ref{prop:existence:pathwisely},
 GFOU $\{Y_t\}_{t\ge0}$ with the initial value $Y_0\in L^1(\Omega)$
 satisfies the stochastic differential equation; 
\begin{eqnarray}
\label{eq:SDE}
dY_t=Y_{t-}dU_t+dB_t^H,
\end{eqnarray}
where $\mathcal{E}(U_t)=e^{-\xi_t}$. 
\end{prop}

\noindent
{\it Proof of Proposition \ref{prop:SDE:GFOU}}. \\
Since the condition of Theorem \ref{thm:rmintegra:p-variation} is satisfied,
the integral $\int_0^tB_{s-}^Hde^{\xi_s}$ also exists in the 
Riemann-Stieltjes path-wise sense. We
use the integration by parts formula to $Y_t$ and obtain
$$
Y_t=e^{-\xi_t}\left(Y_0-\int_0^t B_{s-}^H de^{\xi_s}\right)+B_t^H.
$$ 
If we put $Q_t:=Y_t-B_t^H$, the equation above becomes
$$
Q_t=e^{-\xi_t}\left(Q_0-\int_0^tB_{s-}^H de^{\xi_s}\right), 
$$
where $Q_0=Y_0$. 
Since $\{e^{\xi_t}\}_{t\ge0}$ is a semimartingale and $\{B_t^H\}_{t\ge
0}$ is continuous and adapted, the process $\{Q_t\}_{t\ge0}$ is also
semimartingale. We set $R_t:=e^{-\xi_t}$ and $S_t:=Q_0-\int_0^t
B_{s-}^Hde^{\xi_s}$. Then the integration by parts formula for
semimartingales (e.g. Corollary 2, \Rnum{2} of \cite{Protter:2004}) yields 
\begin{eqnarray*}
Q_t-Q_0 &=& R_tS_t-R_0S_0 \\
        &=& \int_{0+}^t R_{s-}dS_s+\int_{0+}^t S_{s-}dR_s+[R,S]_t-R_0S_0
	 \\
&=& - \int_{0+}^te^{-\xi_{s-}}B_{s-}de^{\xi_s}+\int_{0+}^t
 Q_{s-}e^{\xi_{s-}}de^{-\xi_s}-\int_{0+}^t B_{s-}^H d[e^{-\xi},e^{\xi}]_s \\ 
&=& \int_{0+}^t Q_{s-}dU_s-\int_{0+}^tB_{s-}^H
 \left(e^{-\xi_{s-}}de^{\xi_s}-d[e^{-\xi},e^\xi]_s\right). 
\end{eqnarray*}
Observe the relation between $e^{-\xi_t}$ and $U_t$;
\begin{eqnarray*}
1 &=& e^{\xi_s}e^{-\xi_s} \\
  &=& \int_{0+}^t 
   e^{\xi_{s-}}de^{-\xi_s}+\int_{0+}^te^{-\xi_{s-}}de^{\xi_{s}}+[e^{\xi},e^{-\xi}]_t   \\
  &=& \int_{0+}^t dU_s +
  \int_{0+}^te^{-\xi_{s-}}de^{\xi_{s-}}+[e^\xi,e^{-\xi}]_t. 
\end{eqnarray*}
Using this we obtain 
$$
Q_t-Q_0=\int_{0+}^t\left(Q_{s-}+B_{s-}^H\right)dU_s,
$$
which is equivalent to
$$
Y_t-Y_0=\int_{0+}^t Y_{s-} dU_s+B_t^H.
$$
The proof is now complete. 
\hfill $\Box$ \\

\begin{rem}
The L\'evy measure of $\{\xi_t\}_{t\ge0}$ is obtained from that of $\{U_t\}_{t\ge0}$;
\begin{eqnarray}
\label{eq:relation:levymeasure:xi:U}
\nu_\xi((x,\infty))=\nu_U((-\infty,e^{-x}-1))\quad \text{and}\quad
 \nu_\xi((-\infty,-x))=\nu_U((e^{x}-1,\infty)). 
\end{eqnarray}
 See again Section 2.2 of \cite{erickson:maller:2004}. Hence if $\nu_U$ 
 is concretely given, using criterion of $p$-variation in Section \ref{sec:levy}
 we can check the condition of Proposition \ref{prop:existence:pathwisely}.  
\end{rem}

The following technical Lemma is not difficult but useful for the 
existence of $\{Y_t\}_{t\ge0}$ which is directly constructed from the
stochastic differential equation (\ref{eq:SDE}). The poof is only a
calculation and we omit it.
\begin{lem}
\label{lem:eq:measure:u:xi}
Assume that $\{U_t\}_{t\ge0}$ is a L\'evy process and that
 $\{\xi_t\}_{t\ge0}$ satisfies $\mathcal{E}(U_t)=e^{-\xi_t}$. Then for
 $0<\delta<2$, convergence and divergence of 
 $$
 \int_{|x|<1}|x|^\delta \nu_\xi(dx)\quad \text{and}\quad
 \int_{|x|<1}|x|^\delta \nu_U(dx)
 $$
 are equivalent. 
\end{lem}


\begin{ex}
As an example we consider the stochastic differential equation (\ref{eq:SDE}), where $U_t$ is given by
 an $\alpha$-stable L\'evy motion $\xi_t^\alpha$
 (see Section \ref{sec:levy}). From remark above the L\'evy measure $\nu_\xi$ is given by 
 \begin{equation*}
 \nu_\xi(dx)= \left \{
 \begin{array}{ll}
 c_2\left(1-e^{-x}\right)^{-1-\alpha}e^{-x}dx & \text{on}\quad (0,\infty) \\ 
 c_1\left(e^{-x}-1 \right)^{-1-\alpha}e^{-x}dx & \text{on}\quad (-\infty,0).
 \end{array}
\right.
\end{equation*}
 Observe that
$$
\nu_{\xi}(dx)\sim |x|^{-1-\alpha}dx\quad \text{as}\quad |x|\downarrow 0
$$
and hence variation property of $\xi_t$ is the same as that of
 $U_t=\xi_t^\alpha$. As a result $v_p(\xi)$ is finite if $p>\alpha$. 
\end{ex}

\appendix
\section{Proofs of Section \ref{sec:stationarity:secondorder:GFO}}
\label{sec:pf:secondorder}
{\it Proof of Proposition \ref{thm:second-order-behavior-stationary:GFOU}}. \\
From the stationary version $\overline{Y}$ is definable. By
virtue of the stationarity of $\overline{Y}$ and
Fubini's theorem, we have 
\begin{eqnarray*}
\mathrm{Cov}(\overline{Y}_t,\overline{Y}_{t+s}) &=& \mathrm{Cov}(\overline{Y}_t,\overline{Y}_{t+s}) \\
&=& c_H \int_{-\infty}^0 \int_{-\infty}^0 
 E_\xi[e^{-(\xi_s-\xi_{u-}-\xi_{v-})}]|u-v|^{2H-2}dudv \\
 && +c_H \int_{-\infty}^0 \int_0^s   
E_\xi[e^{-(\xi_s-\xi_{u-}-\xi_{v-})}]|u-v|^{2H-2}dudv \\
&:=& \text{\Rnum{1}}+\text{\Rnum{2}}.
\end{eqnarray*}
First we consider the integral \Rnum{1}. The independent increments
property of $\{\xi_t\}_{t\in\mathbb{R}}$ gives 
\begin{eqnarray*}
E[e^{-(\xi_s-\xi_{u-}-\xi_{v-})}] &=& E[e^{-(\xi_s-\xi_{u}-\xi_{v})}] \\
&=& E[1_{\{u\ge v\}} e^{-(\xi_s-\xi_0)-2(\xi_0-\xi_{u})-(\xi_u-\xi_v)}] \\
&& + E[1_{\{u<v\}} e^{-(\xi_s-\xi_0)-2(\xi_0-\xi_{v})-(\xi_v-\xi_u)}] \\
&=& 1_{\{u\ge v\}}
 E[e^{-(\xi_s-\xi_0)}]E[e^{-2(\xi_{(-u)}-\xi_0)}]E[e^{-(\xi_{u-v}-\xi_0)}]   \\
&&+ 1_{\{u<v\}}E[e^{-(\xi_s-\xi_0)}]
 E[e^{-2(\xi_{(-v)}-\xi_0)}]E[e^{-(\xi_{v-u}-\xi_0)}]. 
\end{eqnarray*}
The integrand $E[e^{-(\xi_s-\xi_{u-}-\xi_{v-})}]$ is symmetric with
respect to $u$ and $v$, and hence
\begin{equation*}
\text{\Rnum{1}} = 2c_H e^{-\theta_1s} \int_{-\infty}^0 \int_{-\infty}^0 
1_{\{u\ge v\}} e^{\theta_2u+\theta_1(v-u)} |u-v|^{2H-2}dudv. 
\end{equation*}
Then further calculation shows
\begin{eqnarray}
\text{\Rnum{1}} &=& 2c_H e^{-\theta_1s} \int_{-\infty}^0 \int_{-\infty}^0 
1_{\{u\ge v\}} e^{\theta_2u+\theta_1(v-u)} |u-v|^{2H-2}dudv \nonumber \\
&& (\text{By change of variables}; x=u-v) \nonumber \\
&=& 2c_H e^{-\theta_1 s} \int_{-\infty}^0 e^{\theta_2u}du \int_0^\infty
 e^{-\theta_1x}x^{2H-2}dx \label{eq:calculation:cov1} \\
&=& \frac{2 c_He^{-\theta_1s}}{\theta_2\theta_1^{2H-1}}\Gamma(2H-1). \nonumber
\end{eqnarray}
Next we consider \Rnum{2}. A similar conclusion as above gives
\begin{eqnarray*}
E[e^{-(\xi_s-\xi_{v-})-(\xi_0-\xi_{u-})}] &=&
 E[e^{-(\xi_s-\xi_v)-(\xi_0-\xi_u)}] \\
&=& e^{-\theta_1(s-v)+\theta_1(u+v)} \\
&=& e^{-\theta_1s+\theta_1(u+v)}
\end{eqnarray*}
and we have 
\begin{eqnarray}
\text{\Rnum{2}} &=& c_He^{-\theta_1 s} \int_{-\infty}^0 e^{\theta_1u} \int_0^s
 (u-v)^{2H-2}dudv \label{eq:calculation:cov2} \\
&& (\text{By change of variables} ;x=v-u ) \nonumber \\
&=& c_H e^{-\theta_1s} \int_{-\infty}^0 e^{\theta_1u}
 \int_{s-u}^{-u}e^{\theta_1(u+x)}x^{2H-2}dudx \\
&=& c_He^{-\theta_1s}\int_0^\infty \int_{-\infty}^0
 1_{\{-x<u<s-x\}}e^{2\theta_1u+\theta_1x}x^{2H-2}dx \nonumber \\
&=& c_H e^{-\theta_1 s}\left(
\int_0^s e^{\theta_1 x}x^{2H-2}dx \int_{-x}^0 e^{2\theta_1u}du 
+
\int_s^\infty e^{\theta_1x}x^{2H-2}dx \int_x^{s-x}e^{2\theta_1u}du
\right) \nonumber \\
&=& c_H e^{-\theta_1s}\bigg\{
\int_0^s\left(\frac{1}{2\theta_1}-\frac{1}{2\theta_1}e^{-2\theta_1x}
	\right)e^{\theta_1x}x^{2H-2}dx \nonumber \\
&&+
\int_s^\infty
\left(\frac{1}{2\theta_1}e^{2\theta_1(s-x)}-\frac{1}{2\theta}e^{-2\theta_1x}\right)
e^{\theta_1x}x^{2H-2}dx 
\bigg\} \nonumber \\
&=&
c_He^{-\theta_1s}\bigg(
\frac{1}{2\theta_1}\int_0^s
e^{\theta_1x}x^{2H-2}dx-\frac{1}{2\theta_1}\int_0^\infty
e^{-\theta_1x}x^{2H-2}dx\\
&& + \frac{e^{2\theta_1 s}}{2\theta_1}\int_s^\infty e^{-\theta_1x}x^{2H-2}dx
\bigg) \nonumber \\
&=& c_He^{-\theta_1s} \bigg( 
\frac{s^{2H+1}}{2\theta_1(2H-1)} {}_1\mathrm{F}_1(2H-1,2H;\theta_1s)
-\frac{1}{2\theta_1^{2H}}\Gamma(2H-1) \nonumber \\
&& + \frac{e^{2\theta_1s}}{2\theta_1^{2H}} \Gamma(2H-1,\theta_1s) 
\bigg) \nonumber 
\end{eqnarray}
Combining \Rnum{1} and \Rnum{2}, we obtain the first assertion.

The series representation of the incomplete Gamma function, i.e.
\begin{eqnarray*}
\Gamma(\alpha,x)=x^{\alpha-1}e^{-x}\left(1+\sum_{n=1}^\infty
\frac{1}{x^n}(\alpha-1)\cdots(\alpha-n)
\right)
\end{eqnarray*}
gives 
\begin{equation}
\label{pf:expansion:gamma}
\Gamma(2H-1,\theta_1s)=\frac{s^{2H-1}}{2\theta_1(H-1)} \sum_{n=0}^\infty 
\prod_{k=1}^{n+1}(2H-k)(\theta_1s)^{-(n+1)} 
\end{equation}
We apply the binomial expansion
\begin{equation*}
(1-u/s)^{2H-2}=\sum_{n=0}^\infty 
\frac{(2H-2)(2H-3)\cdots(2H-n-1)}{n!}\left(-\frac{u}{s}\right)^n, \quad (0<u<s)
\end{equation*}
to the representation
\begin{eqnarray*}
\frac{e^{-\theta_1s}s^{2H-1}}{2\theta_1(H-1)}{}_1\mathrm{F}_1(2H-1,2H;\theta_1s)
&=& \frac{e^{-\theta_1s}}{2\theta_1}\int_0^s e^{\theta_1 x}x^{2H-2}dx \\
&& (\text{By change of variables};u=s-x) \\
&=& \frac{s^{2H-2}}{2\theta_1}\int_0^s e^{-\theta_1u}(1-u/s)^{2H-2}du. 
\end{eqnarray*}
Then we exchange the infinite sum and integral by usual Fubuni's
theorem and obtain
\begin{eqnarray}
&&
\frac{e^{-\theta_1s}s^{2H-1}}{2\theta_1(H-1)}{}_1\mathrm{F}_1(2H-1,2H;\theta_1s)\nonumber \\
&&\quad =-\frac{s^{2H-1}}{2\theta_1(2H-1)}\sum_{n=0}^\infty
\prod_{k=1}^{n+1}(2H-k)(-s\theta_1)^{-(n+1)}\frac{\gamma(n+1,\theta_1s)}{n!}. \label{pf:expansion:f}
\end{eqnarray}
Thus substituting these expansions (\ref{pf:expansion:gamma}) and (\ref{pf:expansion:f}) in
the previous representation of covariance we obtain the result. \hfill
$\Box$ \\

\noindent
{\it Proof of Corollary \ref{cor:second-order-behavior:GFOU}}. \\
Since $X$ is independent it follows that 
\begin{equation}
E[Y_t]E[Y_{t+s}]=(E[X])^2 e^{-\theta_1(2t+s)}. \label{eq:pf:cor:second:1}
\end{equation}
We divide $E[Y_t,Y_{t+s}]$ into piece as follows.
\begin{eqnarray*}
&& E[Y_tY_{t+s}] \\
&& = E\left[
\bigg\{
Xe^{-\xi_t}+\int_0^t e^{-(\xi_t-\xi_{u-})}dB_u^H
\bigg\}
\bigg\{
Xe^{-\xi_{t+s}}+\int_0^{t+s}e^{-(\xi_{t+s}-\xi_{v-})}dB_v^H
\bigg\}
\right] \\
&& = E[X^2 e^{-(\xi_t+\xi_{t+s})}] + E\left[
Xe^{-\xi_{t+s}} + \int_0^{t+s}e^{-(\xi_{t+s}-\xi_{v-})}dB_v^H 
\right] \\
&&\quad +E\left[
Xe^{-\xi_{t+s}}\int_0^t e^{-(\xi_t-\xi_{u-})}dB_u^H
\right]+E\left[
\int_0^t\int_0^{t+s} e^{-(\xi_t-\xi_{u-}+\xi_{t+s}-\xi_{v-})}dB_u^HdB_v^H
\right].
\end{eqnarray*}
The first term is calculated as 
\begin{eqnarray}
E[X^2e^{-(\xi_t+\xi_{t+s})}] &=& E[X^2]E[e^{-(\xi_{t+s}-\xi_t+2\xi_t)}]
\nonumber \\
&=& E[X^2]E[e^{-\xi_s}]E[e^{-2\xi_t}] \nonumber \\
&=& E[X^2]e^{-\theta_1s-\theta_2t}. \label{eq:pf:cor:second:2}
\end{eqnarray}
From Theorem
\ref{thm:second-order-behavior-stationary:GFOU}, the second and the
fourth terms are $0$. We only need the last term, which is calculated as
\allowdisplaybreaks{
\begin{eqnarray}
&& \mathrm{Cov}(\overline{Y}_t,\overline{Y}_{t+s}) - 
E\left[
\int_{-\infty}^0 e^{-(\xi_t-\xi_{u-})}dB_u^H \int_t^{t+s}e^{-(\xi_{t+s}-\xi_{v-})}dB_v^H
\right] \nonumber \\
&& \quad \quad - 
E\left[
\int_{-\infty}^0 e^{-(\xi_t-\xi_{u-})}dB_u^H
\int_{-\infty}^{t+s}e^{-(\xi_{t+s}-\xi_{v-})}dB_v^H  
\right] \nonumber \\
&& \quad \quad -
E\left[
\int_0^t e^{-(\xi_t-\xi_{u-})}dB_u^H \int_{-\infty}^0 e^{-(\xi_{t+s}-\xi_{v-})}dB_v^H
\right] \nonumber \\
&&\quad =
 \mathrm{Cov}(\overline{Y}_t,\overline{Y}_{t+s})-\left\{
e^{-\theta_1t}\mathrm{Cov}(\overline{Y}_0,\overline{Y}_{t+s})
-e^{-\theta_1(t+s)}\mathrm{Cov}(\overline{Y}_0,\overline{Y}_{t})
						 \right\} \nonumber \\
&& \quad \quad -e^{-\theta_1s}E\left[
e^{-\xi_t}\overline{Y}_0\overline{Y}_t
\right] - e^{-\theta_1s}E[e^{-\xi_t}\overline{Y}_0Y_t^{0}],\label{eq:pf:cor:second:3}
\end{eqnarray}
where $Y_t^{0}$ is $Y_t$ with initial value $0$. 

%
By adding up (\ref{eq:pf:cor:second:1}), (\ref{eq:pf:cor:second:2}) and 
(\ref{eq:pf:cor:second:3}), we have
%
\begin{eqnarray*}
&& \mathrm{Cov}(Y_t,Y_{t+s}) \\
&&\ = \mathrm{Cov}(\overline{Y}_t,\overline{Y}_{t+s})-e^{-\theta_1t}
 \mathrm{Cov}(\overline{Y}_0,\overline{Y}_{t+s}) \\
&&\quad\quad +e^{-\theta_1s} \bigg(
E[X^2]e^{-\theta_1t}-(E[X])^2e^{-2\theta_1t}+e^{-\theta_1t}
\mathrm{Cov}(\overline{Y}_0,\overline{Y}_{t}) \\
&& \quad\quad \quad\quad 
-E[e^{-\xi_t}\overline{Y}_0\overline{Y}_t]-E[e^{-\xi_t}\overline{Y}_0Y_t^{0}]
\bigg). 
\end{eqnarray*}
Hence via Theorem \ref{thm:second-order-behavior-stationary:GFOU} the
conclusion follows. \hfill $\Box$ \\

{\bf Acknowledgment:} \ The authors are grateful to Prof.\ Makoto
Maejima for very helpful comments. We give deep thanks to Prof.
Claudia Kl\"uppelberg for careful reading our paper and for useful
comments. Prof.\ Alexander Lindner suggested the topic of Section
\ref{sec:sde:gfou}.  
We also thankful to Prof.\ Jean Jacod for pointing out several
theoretical mistakes and for giving us nice suggestions. 
We shall never forget their kindness.  

\bibliographystyle{plainnat}
\bibliography{references}

\end{document}